\documentclass{amsart}
\usepackage{amssymb}
\usepackage{xcolor}
\usepackage[shortlabels]{enumitem}
\usepackage[citestyle=numeric]{biblatex}
\usepackage{hyperref}
\title[]{Norm-ideal perturbations of one-parameter semigroups and applications}

\author[Boulton]{Lyonell Boulton*} 
\address{*Department of Mathematics and Maxwell Institute for Mathematical Sciences, Heriot-Watt University, Edinburgh, EH14 4AS, UK.}
\email{L.Boulton@hw.ac.uk}

\author[Dimoudis]{Spyridon Dimoudis**}
\address{**School of Mathematics and Statistics, University of St Andrews, Mathematical Institute, North Haugh, St Andrews KY16 9SS, UK.}
\email{sd307@st-andrews.ac.uk}

\date{October 2022}

\newtheorem{Def}{Definition}[section]
\newtheorem{Thm}{Theorem}[section]
\newtheorem{Lma}{Lemma}[section]
\newtheorem{Cor}{Corollary}[section]

\newtheorem{Rem}{Remark}[section]

\renewcommand{\Re}{\operatorname{Re}}

\begin{document}

\begin{abstract}
The notion of equivalence classes of generators of one-parameter semigroups based on the convergence of the Dyson expansion can be traced back to the seminal work of Hille and Phillips, who in Chapter XIII of the 1957 edition of their Functional Analysis monograph, developed the theory in minute detail. Following their approach of regarding the Dyson expansion as a central object, in the first part of this paper we examine a general framework for perturbation of generators relative to the Schatten-von Neumann ideals on Hilbert spaces. This allows us to develop a graded family of equivalence relations on generators, which refine the classical notion and provide stronger-than-expected properties of convergence for the tail of the perturbation series. We then show how this framework realises in the context of non-self-adjoint Schr{\"o}dinger operators.    
\end{abstract}

\maketitle

\tableofcontents

\newpage

\section{Introduction and motivation}

In this work, $\mathcal{H}$ is a separable infinite-dimensional Hilbert space.  For $1\leq q< \infty$, we write $\mathcal{C}_q\equiv \mathcal{C}_{q}(\mathcal{H})$  to denote the $q$-Schatten-von Neumann class of operators and $\mathcal{C}_{\infty}\equiv \mathcal{C}_{\infty}(\mathcal{H})$ to denote the class of compact operators on $\mathcal{H}$. Recall that
$\mathcal{C}_{q}$ and $\mathcal{C}_{\infty}$ are operator ideals in $\mathcal{L}(\mathcal{H})$, the algebra of bounded operators on $\mathcal{H}$, and that
\[
      \mathcal{C}_p\subset \mathcal{C}_q\subset 
     \mathcal{C}_\infty
\]  
for all $p<q$. We will write $\|\cdot\|_q$ for the norm of $\mathcal{C}_{q}$ and unambiguously $\|\cdot\|_{\infty}\equiv \|\cdot\|$.

Let $A$ be the generator\footnote{Here and elsewhere we will say that $A$ is the generator of a one-parameter semigroup, if $-A$ is so in the usual sense \cite[Chapter~X and \S 10.6]{HP1957}.} of a $C_0$ one-parameter semigroup $T(t,A)=\mathrm{e}^{-At}$ on $\mathcal{H}$. Classically \cite{HP1957}, we know conditions for an operator $B:\operatorname{Dom}(A)\longrightarrow \mathcal{H}$ to ensure that $A+B$ is also the generator of a $C_0$ one-parameter semigroup and that $T(t,A+B)$ is given in terms of $T(t,A)$ by a  (Dyson) expansion convergent in the operator norm. One of our main goals below is to find additional conditions on $B$, so that $T(t,A+B)-T(t,A)\in \mathcal{C}_q$ and that the Dyson formula has a tail convergent in $\mathcal{C}_r$, for suitable $r$. Here $q$ and $r$ are related but can be different. Our results significantly improve those of the work \cite{B2019}. 

In order to avoid technical difficulties with integrability of families of operators, hence possibly distracting from our main purpose, we focus  on immediately norm continuous semigroups. This does not compromise the general nature of our results, as they cover important cases, such as those of Gibbs type \cite{Z2019} as well as those satisfying the property that $T(t,A)f\in \operatorname{Dom}(A)$ for all $f\in \mathcal{H}$ and $t>0$ \cite[Theorem~1.28]{D1980}. This includes those that can be extended to bounded holomorphic semigroups on a sector \cite[Theorem~2.39]{D1980}. As a proof of concept of the framework, in the second part of the paper we show applications of the abstract setting to heat semigroups associated with non-self-adjoint Schr{\"o}dinger operators.

\section{Perturbation of generators relative to an operator ideal}

We begin with the classical definitions found in \cite[Chapter~XIII]{HP1957} re-written in contemporary language. Let $A$ be the generator of a $C_0$ one-parameter semigroup $T(t,A)=\operatorname{e}^{-At}$. For a linear operator $B:\operatorname{Dom}(A)\longrightarrow \mathcal{H}$ we say  $B\in \mathcal{J}(A)$ iff $B$ is relatively bounded with respect to $A$. Recall that this is equivalent to having $BR(\lambda,A)\in \mathcal{L}(\mathcal{H})$ for some (and hence all) $\lambda\not \in \operatorname{Spec}(A)$. Here and everywhere below
\[
    R(\lambda,A)=(\lambda-A)^{-1}
\] 
is the resolvent operator.

For $B\in \mathcal{J}(A)$ there exists a unique extension of $B$,  \cite[Theorem~13.3.1]{HP1957} denoted below by $\tilde{B}$, such that 
\[x\in \operatorname{Dom}(\tilde{B}) \qquad \iff \qquad  y:=\lim_{\lambda \rightarrow -\infty}|\lambda|BR(\lambda,A)x  \text{ exists},\]
and in this case $y=\tilde{B}x$. By virtue of the Hille-Yosida Theorem, $\|R(\lambda,A)\|$ decays linearly as $\lambda\to -\infty$, therefore, indeed \[\operatorname{Dom}(B)=\operatorname{Dom}(A)\subset \operatorname{Dom}(\tilde{B
}).\] 
The extension $\tilde{B}$ is also relatively bounded with respect to $A$, it has the same relative bound as that of $B$ and it plays a significant role in the theory of perturbations of one-parameter semigroups, as we shall see next. 

Two well-known observations are now in place. Firstly, if $B\in \mathcal{J}(A)$, the operator
\[BT(t,A)\!\upharpoonright _{\operatorname{Dom}(A)}:\operatorname{Dom}(A)\longrightarrow \mathcal{H}\]
is always well defined by elementary properties of $C_0$-semigroups \cite[Lemma~1.1]{D1980}. If $BT(t,A)\!\upharpoonright _{\operatorname{Dom}(A)}$ is bounded taking the operator norm in $\operatorname{Dom}(A)$, then \[\tilde{B}T(t,A):\mathcal{H}\longrightarrow \mathcal{H}\] is well defined, bounded and
\[ \|  \tilde{B}T(t,A)\|=\|BT(t,A)\!\upharpoonright_{\operatorname{Dom}(A)}\|.\] 
Secondly, the operator $B\in \mathcal{J}(A)$ may or may not be closable, but if it is closable, then  \[ \operatorname{Dom}(B)=\operatorname{Dom}(A)\subseteq \operatorname{Dom}(\tilde{B})\subseteq
\operatorname{Dom}(\overline{B}).\] 

The Hille-Phillips class $\mathcal{B}(A)$ (or class $\mathcal{P}$ as referred to in \cite{D1980} for closable perturbations) is defined as follows \cite[Definition~13.3.5]{HP1957}. A linear operator $B\in\mathcal{B}(A)$ iff
\begin{enumerate}[i)]
\item \label{ci} $B \in \mathcal{J}(A)$,
\item \label{cii} $BT(t,A)\!\upharpoonright _{\operatorname{Dom}(A)}$ is bounded for all $t>0$,
\item \label{inteforB(A)} $\int_{0}^{1} \|\tilde{B}T(s,A)\| \,\mathrm{d}s< \infty$.
\end{enumerate}
The conditions \ref{ci} and \ref{cii} imply that the map $t\longmapsto \tilde{B}T(t,A)$ is continuous in the strong operator topology \cite[Lemma~13.3.3]{HP1957}. Thus, the scalar map $t\longmapsto \|\tilde{B}T(t,A)\|$ is Borel measurable and the integral \ref{inteforB(A)} is well defined. Note that the class $\mathcal{B}(A)$ is a linear subspace of $\mathcal{J}(A)$.

 If $B\in \mathcal{B}(A)$, then \[
A+B:\operatorname{Dom}(A)\longrightarrow \mathcal{H} 
\]
is the generator of a $C_0$ one-parameter semigroup. In \cite{B2019} we determined additional conditions on a closable operator $B\in \mathcal{B}(A)$ so that, when $A$ is the generator of a (Gibbs) semigroup $T(t,A)\in\mathcal{C}_1$ for all $t>0$, then also $T(t,A+B)\in\mathcal{C}_1$ for all $t>0$. A key ingredient to ensure the latter \cite[Lemma~1]{B2019} is that the integral in \ref{inteforB(A)} converges in the norm of $\mathcal{C}_q$ for some $q<\infty$, so that the Dyson expansion for $T(t,A+B)$ is not only convergent in the operator norm, but also in the stronger norm of the $q$-Schatten-von Neumann class. The definition after the next lemma, and the subsequent results, extend this observation and the findings of \cite[Section~2]{B2019} to a significant level of generality.
 
A $C_0$ one-parameter semigroup is often called immediately norm continuous, if $t\mapsto T(t,A)$ is continuous in the operator norm for all $t>0$. We will adhere to this terminology and will consider from now on that $T(t,A)$ is an immediately norm continuous semigroup.

\begin{Lma} \label{lemma1}
Let $B\in \mathcal{J}(A)$ be such that $\tilde{B}T(t,A)\in \mathcal{C}_{q}$ for all $t>0$. If $T(t,A)$ is immediately norm continuous, then $\tilde{B}T(t,A)$ is a continuous function with respect to  $\|\cdot\|_q$ for all $t>0$.
\end{Lma}
\begin{proof}
Set $t$ and $t_0$, such that either $t\geq t_{0}>0$ or $t_{0}>t>t_{0}/2>0$.  Then
\begin{align*}
\|\tilde{B}T(t,A)-\tilde{B}T(t_{0},A)\|_{q}&=\|\tilde{B}T(t_{0}/2,A)(T(t-t_{0}/2,A)-T(t_{0}/2,A))\|_{q} \\ & \leq\|\tilde{B}T(t_{0}/2,A)\|_{q}\|T(t-t_{0}/2,A)-T(t_{0}/2,A)\|_{\infty}.
\end{align*}
By taking the limit $t \rightarrow t_{0}$ and applying the hypothesis at $t_{0}/2$, the claim follows. 
\end{proof}

Under the hypotheses of this lemma, $\tilde{B}T(t,A)$ is almost separably-valued. Indeed, the image of the separable set $(0,\infty)$ under this continuous function is separable. Hence, by virtue of the Pettis Theorem, this family of operators is strongly measurable on $(0,\infty)$.

\begin{Def} \label{d3}
Let $1\leq q< \infty$ and let $A$ be the generator of an immediately norm continuous semigroup $T(t,A)$. We will write that a linear operator $B\in\mathcal{B}_{q}(A)$ iff 
\begin{enumerate}[a)]
\item \label{a)forBq(A)}$B\in \mathcal{B}(A)$,
\item \label{b)forBq(A)} $\tilde{B}T(t,A)\in \mathcal{C}_{q}$ for all $t>0$,
\item \label{c)forBq(A)} $\int_{0}^{1} \|\tilde{B}T(s,A)\|_{q}\,\mathrm{d}s< \infty$.
\end{enumerate}
\end{Def}

If $T(t,A)$ is immediately norm continuous and $B\in \mathcal{B}(A)$, then  $T(t,A+B)$ is also immediately norm continuous \cite[Theorem~13.4.2]{HP1957}. 
Hence, perturbations by the class $\mathcal{B}_q(A)$ also preserve this property. Below we will often use this fact without further mention. 

Note that $\mathcal{B}_{q}(A)$ is a linear subspace of $\mathcal{J}(A)$. Indeed, let $B,C \in \mathcal{B}_{q}(A)$. 
Then $(B+C)\in\mathcal{B}(A)$ by linearity of $\mathcal{B}(A)$. Moreover
\[
       \widetilde{(B+C)}T(t,A)x=(\tilde{B}+\tilde{C})T(t,A)x \qquad \forall x\in \mathcal{H},
\]
because the equality holds for all $x\in\operatorname{Dom}(A)$ and the latter is dense in $\mathcal{H}$. Thus
\[\widetilde{(B+C)}T(t,A)=\tilde{B}T(t,A)+\tilde{C}T(t,A)\in \mathcal{C}_q,\] ensuring \ref{b)forBq(A)}. Hence
\[
    \|\widetilde{(B+C)}T(t,A)\|_q\leq\|\tilde{B}T(t,A)\|_q+\|\tilde{C}T(t,A)\|_q
\]
ensures \ref{c)forBq(A)} for the sum. Consequently, $(B+C)\in \mathcal{B}_q(A)$. The fact that multiplying by a scalar preserves membership to $\mathcal{B}_q(A)$ is obvious.

The upper limit in the integration \ref{c)forBq(A)} could be set to any other $t>1$ without altering the class identified in Definition~\ref{d3}. Indeed, from the semigroup property, it follows that
\begin{align*}
 \int_{n-1}^{n} \|\tilde{B}T(s,A)\|_{q}\,\mathrm{d}s&= \int_{0}^{1} \|\tilde{B}T(s,A)T(n-1,A)\|_{q}\,\mathrm{d}s  \\ &\leq\|T(n-1,A)\|_{\infty}\int_{0}^{1} \|\tilde{B}T(s,A)\|_{q}\,\mathrm{d}s,
\end{align*}
for all $n\in \mathbb{N}$. 
Therefore 
\begin{align*}
 \int_{0}^{t} \|\tilde{B}T(s,A)\|_{q}\,\textnormal{d}s<\infty,
\end{align*}
for any $t>1$ if and only if \ref{c)forBq(A)} holds true. Moreover, if the semigroup has a negative growth bound, then we can take $t=\infty$.

It is classically known that $\mathcal{B}(A)$ yields a partition of the family of generators of $C_0$ semigroups on a given Hilbert space. We will show below that the $\mathcal{B}_q(A)$ families also determine nested partitions of the set of generators of immediately norm continuous semigroups. These partitions correspond to the equivalence classes associated to the following relations.
\begin{Def} \label{equivalentclass}
Let $1\leq q< \infty$. Let $A_{1}$ and $A_{2}$ be generators of immediately norm continuous semigroups. We will write $A_{1} \approx_q A_{2}$ iff $A_{2}=A_{1}+B$ for $B \in \mathcal{B}_{q}(A_{1})$.
\end{Def}

The proof that these relations are equivalences will be given in Section~\ref{proofofeq}. 

\section{Perturbation formulas for semigroups}

In this section we derive the convergence in norm $\|\cdot\|_{q}$ of the classical perturbation identities for semigroups. One of our main tools is the next lemma about continuity of convolutions in suitable norms, which might be well known. As it will be crucial in various arguments later on, we include a self-contained proof in the context of the trace ideals. 

Firstly recall the fundamental interpolation identity \cite[Lem.XI.9.20]{DS1963}. For $S\in \mathcal{C}_q$ and $T\in \mathcal{C}_r$,
\begin{equation}  \label{Holder}
   \|ST\|_p\leq\|S\|_q\|T\|_r      ,
\end{equation}
whenever $p,\,q,\,r\geq 1$ are such that $\frac{1}{p}=\frac{1}{q}+\frac{1}{r}$.

\begin{Lma} \label{lem2}
Let $1\leq p,\,q,\,r\leq \infty$ be such that $\frac{1}{p}=\frac{1}{q}+\frac{1}{r}$. Let $F,G: (0,\infty) \longrightarrow \mathcal{L}(\mathcal{H})$ be operator-valued functions such that $F$ is continuous in $\|\cdot\|_{q}$ and $G$ is continuous in $\|\cdot\|_{r}$. Suppose that \[
\int_{0}^{t} \|F(s)\|_{q} \,\textnormal{d} s<\infty \qquad \text{and} \qquad \int_{0}^{t} \|G(s)\|_{r}\, \textnormal{d} s<\infty,
\] for all $t>0$. Then
the function $s\mapsto F(t-s)G(s)$ is continuous in $\|\cdot\|_{p}$ for all $0<s<t$ and \[\int_{0}^{t} \|F(t-s)G(s)\|_{p}\,\textnormal{d}s < \infty.\] 
Set
\[
     H(t):=(F*G)(t)=\int_{0}^{t} F(t-s)G(s) \,\textnormal{d}s.
\]
If $p<\infty$, then $H:(0,\infty) \longrightarrow \mathcal{C}_{p}$ and it is continuous in the norm $\|\cdot\|_{p}$. 
\end{Lma}
\begin{proof}
Without loss of generality assume $1\leq q \leq r \leq \infty$.
Fix $t>0$ and let $0<s,s_{0}<t$. By the triangle inequality and \eqref{Holder}, we get 
\begin{align*}
\|F(t-s)&G(s)-F(t-s_{0})G(s_{0})\|_{p} \\& \leq 
\|F(t-s)(G(s)-G(s_0))\|_p+\|(F(t-s)-F(t-s_0))G(s_0)\|_{p} \\ &\leq
\|F(t-s)\|_{q}\|G(s)-G(s_{0})\|_{r}+ \|F(t-s)-F(t-s_{0})\|_{q}\|G(s_0)\|_{r}.
\end{align*}
This, and the continuity of $F$ and $G$, ensure the first conclusion. In particular, $\|F(t-s)G(s)\|_{p}$ is measurable.

Let us turn our attention to the second claim. For $a,b>0$, let
\begin{equation*}
f(a,b):=\max_{a\leq x \leq b} \|F(x)\|_{q} \quad \text{and} \quad  
g(a,b):=\max_{a\leq x \leq b} \|G(x)\|_{r}.
\end{equation*}
From the continuity of $F$ and $G$, it follows that both these quantities are finite and continuous in the two variables.
Hence, writing the integral into two parts and applying \eqref{Holder}, yields
\begin{align*}
\int_{0}^{t}\|F(t-s)G(s)\|_{p}\,\mathrm{d}s &\leq f(t/2,t)
\int_{0}^{t/2}\|G(s)\|_{r}\,\mathrm{d}s  \\  &+g(t/2,t)\int_{t/2}^{t}\|F(t-s)\|_{q} \,\mathrm{d}s<\infty,
\end{align*}
as needed.

Now assume $p<\infty$ and consider $H(t)$ as in the hypothesis. From the first two conclusions, it follows that the integral appearing on the right hand side of its definition converges in $\mathcal{C}_{p}$. Hence $H(t)\in \mathcal{C}_{p}$ and it is only left to prove its continuity in $\|\cdot\|_p$.

For that purpose, let $\xi_{0}>0$ and choose $\delta>0$ such that $\xi_{0}-2\delta>0$. Let $\xi_{1}$ be such that $|\xi_{1}-\xi_{0}|\leq \delta$. Separating the integral in an appropriate manner gives
\begin{align*}
\|H(\xi_{1})-H(\xi_{0})\|_{p}&\leq  \int_{0}^{\xi_{0}-2\delta} \big\|(F(\xi_{1}-s)-F(\xi_{0}-s))G(s)\big\|_{p}\, \textnormal{d}s\\ &+ \sum_{j=0}^{1}\int_{\xi_{0}-2\delta}^{\xi_{j}} \|F(\xi_{j}-s)G(s)\|_{p}\, \textnormal{d}s.
\end{align*}
Now, for $0<s<\xi_{0}-2\delta$, 
\begin{equation*}
\big\|(F(\xi_{1}-s)-F(\xi_{0}-s))G(s)\big\|_{p}\leq 2f(\delta,\xi_{0}+\delta) \|G(s)\|_{r}.
\end{equation*}
This ensures that the first integrand, which is continuous in $\xi_{1}$, is dominated by an integrable function. Hence, by the Dominated Convergence Theorem, it converges to zero as $\xi_{1} \rightarrow \xi_{0}$. 
Moreover, for $\xi_{0}-2\delta<s<\xi_{j}$, $j=0,1$,
\begin{equation*}
\|F(\xi_{j}-s)G(s)\|_{p} \leq g(\xi_{0}-2\delta, \xi_{0}+\delta) \|F(\xi_{j}-s)\|_{q}.
\end{equation*}
Then,
\begin{align*}
\lim_{\delta \to 0}\sum_{j=0}^{1} &\int_{\xi_{0}-2\delta}^{\xi_{j}} \|F(\xi_{j}-s)G(s)\|_{p}\, \textnormal{d}s \\ & \leq 2 \lim_{\delta \to 0}\left(g(\xi_{0}-2\delta, \xi_{0}+\delta)  \int_{0}^{3\delta} \|F(s)\|_{q}\, \textnormal{d} s\right)=0.
\end{align*}
From this, the last conclusion follows.
\end{proof}
The next observation about the norm of $H(t)$, which is a consequence of \eqref{Holder}, will be useful below. In the context of the lemma above, let
\[
     \phi(t)=\|F(t)\|_q \qquad \text{and} \qquad \psi(t)=\|G(t)\|_r.
\]
Then 
\begin{equation} \label{convbound}
    \|H(t)\|_p\leq (\phi*\psi)(t),
\end{equation}
irrespective of the fact that the right hand side can be infinity.

Below we will be mainly concerned with the choice
\[
    \phi(t)=\|T(t,A)\|_{\infty} \qquad \text{and} \qquad \psi(t)=\|\tilde{B}T(t,A)\|_q,
\]
for $B\in \mathcal{B}_q(A)$. From the semigroup property, we know that $\phi(t)$ is submultiplicative, that is
$
  \phi(t+s) \leq \phi(t)\phi(s).
$
Also, the growth bound of $T(t,A)$ is
\[
    \omega_0(A)=\lim_{t\to \infty} \frac{\log \phi(t)}{t}.
\]
Moreover, 
\begin{equation} \label{weirdsubad}
\psi(t+s)=\|\tilde{B}T(t+s,A)\|_q \leq \|\tilde{B}T(t,A)\|_q \|T(s,A)\|_\infty=\psi(t)\phi(s).
\end{equation}
Then
\[
    \limsup_{t\to\infty} \frac{\log \psi(t)}{t}\leq \omega_0(A).
\]
Now, part~\ref{c)forBq(A)} of Definition~\ref{d3} can be recast as  \begin{equation} \label{intdyson2} \int_{0}^{1}(\phi(s)+\psi(s))\,\textnormal{d}s<\infty.\end{equation}
In fact, by the classical statement \cite[Lemma~13.4.1]{HP1957}, we know that for all $\omega>\omega_0(A)$,
\begin{equation} \label{expoestimate}
    \int_{0}^\infty \mathrm{e}^{-\omega s} (\phi(s)+\psi(s))\,\textnormal{d}s<\infty.
\end{equation}
If we denote the left hand side of this by $M_{\omega}$, then
\[
     \psi(t)\leq \frac{M_{\omega}^2 \mathrm{e}^{\omega t}}{t^2} \qquad \forall t>0.
\]
These observations will shortly be crucial. 

Let us determine precise convergence properties in a norm $\|\cdot\|_{r}$ of the Dyson expansion of $T(t,A+B)$ whenever $B\in\mathcal{B}_{q}(A)$. Here $r$ and $q$ will be related, but they might not be equal. Our method of proof follows closely those of \cite[Theorem~13.4.1 and its Corollary~1]{HP1957}. Although the mechanism might not be particularly novel, the conclusions are surprising.

We begin by recalling the classical statement. If $B\in\mathcal{B}(A)$, then for $x\in \mathcal{H}$,
\begin{equation} \label{eqDyson}
T(t,A+B)x-T(t,A)x=\sum_{n=1}^{\infty} S_{n}(t)x,
\end{equation}
where 
\begin{equation} \label{multifoldints} 
\begin{aligned}
S_{1}(t)x&=\int_{0}^{t} T(t-s,A)\tilde{B}T(s,A)x\,\textnormal{d}s \\
S_{n}(t)x&=\int_{0}^{t} T(t-s,A)\tilde{B}S_{n-1}(s)x\,\textnormal{d}s \qquad n> 1.
\end{aligned}
\end{equation} 
The series on the right hand side converges in the operator norm, uniformly on $(0,\alpha)$ for all $\alpha>0$. The next two theorems show that the convergence of \eqref{eqDyson} improves gradually in the tail when $B\in \mathcal{B}_q(A)$. Note that, for the case of Gibbs semigroups, these two results are sharper than those formulated in \cite{B2019}.

\begin{Thm} \label{dyson1}
Let $r\geq 1$. Let $T(t,A)$ be an immediately norm continuous one-parameter semigroup. If $B \in \mathcal{B}_{r}(A)$, then the integrals in \eqref{multifoldints} converge
\begin{itemize}
\item in the norm $\|\cdot\|_{r/n}$, so $S_{n}(t) \in \mathcal{C}_{r/n}$  for $n\leq r$, 
\item in the norm $\|\cdot\|_{1}$, so $S_{n}(t) \in \mathcal{C}_{1}$ for $n>r$.\end{itemize} Moreover, these operator functions are continuous in the respective norms.
\end{Thm}
\begin{proof} 
Firstly, consider the case $n=1$. Take $F(t)=T(t,A)\in \mathcal{L}(\mathcal{H})$ and $G(t)=\tilde{B}T(t,A)\in \mathcal{C}_r$ in Lemma~\ref{lem2} with $q=\infty$. Then $S_1(t)\in \mathcal{C}_p$ with $p=r$. This is the claim for $n=1$. Moreover,  $S_1:(0,\infty)\longrightarrow \mathcal{C}_r$ is continuous. 

Now consider the case $n=2$. Since $B\in\mathcal{B}(A)$, then by \cite[Lemma~13.3.5]{HP1957} 
\[
   \tilde{B}S_1(t)x=\int_0^t \tilde{B}T(t-s,A)\tilde{B}T(s,A)x\,\mathrm{d}s,
\]
for all $x\in \mathcal{H}$. Two possibilities arise. 

One is that $r\geq 2$. In this case, take $F(t)=G(t)=\tilde{B}T(t,A)\in \mathcal{C}_{r}$ with $q=r$ in Lemma~\ref{lem2}. Then 
$\tilde{B}S_1(t)\in \mathcal{C}_{r/2}$ is continuous for all $t>0$. Hence, take $F(t)=T(t,A)\in \mathcal{L}(\mathcal{H})$ and $G(t)=\tilde{B}S_1(t)$ with $r/2$ for $G$ in  Lemma~\ref{lem2}, to obtain that $S_2(t)\in \mathcal{C}_{p}$ with $p=r/2$ and that $S_2:(0,\infty)\longrightarrow \mathcal{C}_{r/2}$ is continuous. 

The other possibility is $r<2$. Such being the case, take $F(t)=G(t)=\tilde{B}T(t,A)\in \mathcal{C}_{2}$, which satisfies the hypotheses of Lemma~\ref{lem2} with $q=r=2$, by monotonicity of the Schatten norm. Then, $\tilde{B}S_1(t)\in \mathcal{C}_{1}$. Once again applying Lemma~\ref{lem2} as before gives that $S_2(t)\in \mathcal{C}_{p}$ with $p=1$, and that $S_2:(0,\infty)\longrightarrow \mathcal{C}_{1}$ is continuous. This completes the proof of the case $n=2$. 

For $n\geq 3$, we can proceed in a similar way, showing that $S_{n}(t)\in \mathcal{C}_{r/n}$ for $n\leq r$ or $S_{n}(t)\in \mathcal{C}_{1}$ otherwise, with continuity in the respective norms. The proof can be completed proceeding inductively. We omit further details.
\end{proof}

For $k\in\mathbb{N}$, we denote by
\begin{equation} \label{tailDyson}
      V_k(t)x=\sum_{n=k+1}^\infty S_n(t)x    
\end{equation}
the tail of the series on the right hand side of \eqref{eqDyson}. We now provide precise details about the convergence of \eqref{tailDyson}.

\begin{Thm} \label{dyson2}
Let $q\geq 1$. Let $T(t,A)$ be an immediately norm continuous one-parameter semigroup and let $B \in \mathcal{B}_{q}(A)$. Let $r$ be the integer part of $q$. Then, $V_{r}(t)\in\mathcal{C}_1$ for all $t>0$. Moreover, for all $\alpha>0$, the operator series on the right hand side of \eqref{tailDyson} for $k=r$ converges absolutely in the norm $\|\cdot\|_{1}$ uniformly for $t\in (0,\alpha)$.
\end{Thm}    
\begin{proof}
Let $t>0$. Set $r$ to be the integer part of $q$, as in the hypothesis. Our first goal is to show that
\begin{equation}   \label{tail_trace_convergence}
     \sum_{n=r+1}^\infty \|S_n(t)\|_1<\infty.
\end{equation} 
According to Theorem~\ref{dyson1}, we know that each one of the terms of this series is finite. We aim at applying the classical result \cite[Lemma~13.4.3]{HP1957} as follows.
Set scalar functions $\phi,\psi_k:(0,\infty)\longrightarrow (0,\infty)$ defined by
\[\phi(t)=\|T(t,A)\|_{\infty} \qquad \text{and}\qquad  \psi_k(t)=\|\tilde{B}S_{k-1}(t)\|_{\max\{1,\frac{q}{k}\}} \qquad k\in\mathbb{N},
\]
with the usual convention $S_0(t)=T(t,A)$. The arguments in the proof of Theorem~\ref{dyson1} show that $\psi_k(t)<\infty$ for all $t>0$. Moreover, the hypotheses and Lemma~\ref{lemma1} ensure that these real-valued functions are continuous.

Next, let $n\geq r+1$.  By virtue of \eqref{convbound}, taking $\psi(t)=\psi_n(t)$, it follows that
\[
     \|S_n(t)\|_1 \leq \phi * \psi_n(t).
\]
Also from \eqref{convbound}, it follows that
\[
    \psi_n(t)\leq \psi_1*\psi_{n-1}(t).
\]
Then, as all these functions are positive and convolutions preserve inequalities,
\begin{equation} \label{bound_for_convergence_convolution}
     \|S_n(t)\|_1\leq \phi *\psi_1*\psi_{n-1}(t)\leq \cdots \leq \phi*\psi_1^{[n*]}(t).
\end{equation}
Here $\psi_1^{[n*]}=\psi_1*\cdots*\psi_1$ with the total number of convolutions on the right hand side being $n-1$.

From \eqref{weirdsubad} and \eqref{intdyson2}, it follows by \cite[Lemma~13.4.3]{HP1957}, that the series
\begin{equation}  \label{seriestheta}
    \theta(t)=\sum_{n=r+1}^\infty \phi*\psi_1^{[n*]}(t)<\infty,
\end{equation}
and so \eqref{bound_for_convergence_convolution} implies  \eqref{tail_trace_convergence}. This shows that $V_r(t)\in \mathcal{C}_1$ for all $t>0$. Moreover, the same lemma yields that for all $0<\varepsilon<1$ the convergence is uniform for $t\in (\varepsilon,1/\varepsilon)$. 

We finally show that the convergence is uniform also in $(0,\alpha)$.
Set
\[
    \theta_0(t)=\sum_{n=r+1}^\infty \psi_1^{[n*]}(t),
\]
so that $\theta(t)=\phi*\theta_0(t)$. From \eqref{expoestimate} it follows that for large enough $\omega>\omega_0(A)$, 
\[
    \int_0^\infty \mathrm{e}^{-\omega s} \psi_1(s) \textnormal{d} s<1.
\]
Since 
\[
    \int_0^\infty \mathrm{e}^{-\omega s} \psi^{[n*]}_1(s)\, \textnormal{d}s=
\left(\int_0^\infty \mathrm{e}^{-\omega s} \psi_1(s) \,\textnormal{d} s\right)^n,
\]
we have that
\[
     \int_0^\infty \mathrm{e}^{-\omega s} \theta_0(s) \,\textnormal{d}s<\infty.
\]
This implies that $\int_0^\alpha \theta_0(s)\textnormal{d}s<\infty$ for all $\alpha>0$. Let $M_\alpha>0$ be such that $\phi(t)=\|T(t,A)\|<M_\alpha$ for all $t\in(0,\alpha)$. Then
\[
    \theta(t)=\int_0^t \phi(t-s)\theta_0(s)\textnormal{d}s\leq M_\alpha\int_0^t
\theta_0(s)\textnormal{d}s<\infty,
\]
so that $\lim_{t\to 0}\theta(t)=0$. From the latter, it follows that the series in \eqref{seriestheta} and thus in \eqref{tail_trace_convergence} converge uniformly for all $t\in(0,\alpha)$ as claimed.
\end{proof} 

We conclude this section by highlighting three corollaries, consequence of  Theorem~\ref{dyson2}. The first one is relevant in the context of Gibbs semigroups. It extends \cite[Lemma~1]{B2019} in that $B$ is not required to be a closable operator.

\begin{Cor} \label{gibbs}
Let $T(t,A) \in \mathcal{C}_{1}$ for all $t>0$ and $B \in \mathcal{B}_{q}(A)$. Then $T(t,A+B) \in \mathcal{C}_{1}$ for all $t>0$.
\end{Cor}
\begin{proof}
Since $T(t,A) \in \mathcal{C}_{1}$, $T(t,A)$ is compact for all $t>0$, and by \cite[Theorem 1.30]{D1980} it is immediately norm continuous. The previous theorem applies, and $T(t,A+B)=(T(t,A+B)-T(t,A))+T(t,A) \in \mathcal{C}_{1}$ for all $t>0$. Here we have used the fact that a one-parameter semigroup is in $\mathcal{C}_{1}$ for all $t>0$ if and only if it is in $\mathcal{C}_{q}$ for all $t>0$ and some $1<q<\infty$.
\end{proof}

The second corollary involves the asymptotic behaviour of the perturbed semigroup at the origin, in relation to the unperturbed semigroup.

\begin{Cor} \label{orderofDyson}
Let $T(t,A)$ be an immediately norm continuous one-parameter semigroup. If $B \in \mathcal{B}_{q}(A)$, then there exists a constant $\omega(A,B)\geq \omega_0(A)$ such that
for $\omega>\omega(A,B)$,
\[
     \|T(t,A+B)-T(t,A)\|_q\leq \frac{\mathrm{e}^{\omega t}}{t^2} \qquad \forall t>0.
\] 
\end{Cor}
\begin{proof}
As a consequence of classical estimates for convolutions, \cite[Lemma~13.4.2 and (13.4.7)]{HP1957}, there exists a constant, that we set to be $\omega(A,B)> \omega_0(A)$, such that for $\omega\geq\omega(A,B)$,
\begin{equation} \label{uniformboundS}
      \|S_n(t)\|_{\max\{1,\frac{q}{n}\}}\leq \frac{\mathrm{e}^{\omega t}}{2^nt^2} \qquad \forall t>0.
\end{equation}
From this, it is straightforward to obtain the claim of this corollary.
\end{proof} 

We finally consider Duhamel's formula.  Recall that if $B\in \mathcal{B}(A)$ and $x \in \mathcal{H}$, then
\begin{equation}  \label{Duhamelproper}
   T(t,A+B)x-T(t,A)x=\int_{0}^{t}T(t-s,A+B)\tilde{B}T(s,A)x\, \textnormal{d}s.
\end{equation}

\begin{Cor}  \label{Duhamel}
Let $T(t,A)$ be an immediately norm continuous one-parameter semigroup. If $B\in \mathcal{B}_q(A)$, then the right hand side of \eqref{Duhamelproper} converges in $\|\cdot\|_q$.
\end{Cor}
\begin{proof}
 In Lemma~\ref{lem2} take $F(t)=T(t,A+B)\in \mathcal{L}(\mathcal{H})$ and $G(t)=\tilde{B}T(t,A)\in \mathcal{C}_q$. 
\end{proof}

\section{The equivalence relations} \label{proofofeq}

We are now in the position to show that $\approx_q$ from Definition~\ref{equivalentclass} is an equivalence relation on the class of generators of immediately norm continuous semigroups.

Our first main goal is to determine a generalisation of Theorem~\ref{dyson2}, which is quite interesting in its own right. It relates two perturbations of $A$ which are in different classes.

\begin{Thm} \label{BtildeDyson}
Let $T(t,A)$ be an immediately norm continuous one-parameter semigroup. For $p\geq q\geq 1$, let $B_0\in \mathcal{B}_p(A)$ and $B \in \mathcal{B}_{q}(A)$. Then the following is true. \begin{enumerate}[i)] 
\item \label{BtildeDyson_i)} $\widetilde{B_0} T(t,A+B)\in\mathcal{C}_p$. 
\item \label{BtildeDyson_ii)} To be precise, if $\ell$ is the integer part of $q-\frac{q}{p}$,
\begin{equation}  \label{BtildeDysonformula}
\tilde{B}_0T(t,A+B)=\widetilde{B_0}S_0(t)+\cdots +\widetilde{B_0}S_{\ell}(t)+W(t),
\end{equation}
where 
\[
    \widetilde{B_0}S_n(t)\in \mathcal{C}_{r_n}
\quad \text{and} \quad
    r_n=\frac{pq}{np+q}, \quad n=0, 1,\ldots,l.
\]
\item \label{BtildeDyson_iii)} Moreover,
\[
    W(t)=\sum_{n=\ell+1}^\infty \widetilde{B_0}S_n(t)\in\mathcal{C}_1,
\]     
where the series converges in the norm $\|\cdot\|_1$ uniformly in $t\in (0,\alpha)$ for all $\alpha>0$. 
\item \label{BtildeDyson_iv)} Furthermore,
\[
    \int_0^1\|\tilde{B}_0S_n(s)\|_{r_n}\, \mathrm{d}s<\infty, \qquad 
n=0, 1,\ldots,\ell,
\]
\item \label{BtildeDyson_v)} and
\[
      \int_0^1 \|W(s)\|_1\,\mathrm{d}s<\infty.
\]
\end{enumerate}
\end{Thm}
\begin{proof}
Firstly, note that \ref{BtildeDyson_i)}, follows from \ref{BtildeDyson_ii)} and \ref{BtildeDyson_iii)}.

Let us show \ref{BtildeDyson_ii)}. By virtue of \cite[Lemma~13.5.1]{HP1957}, we know from the fact that $B_0,\,B\in \mathcal{B}(A)$, that for all $x\in \mathcal{H}$
\[
    \widetilde{B_0}S_{n}(t)x=\int_0^t \widetilde{B_0}T(t-s,A)\tilde{B}S_{n-1}(t)x\,\mathrm{d}s,
\]
with the integral converging in operator norm.
We aim at applying Lemma~\ref{lem2} recursively. For $n=0$, note that $\widetilde{B_0} S_0(t)=\widetilde{B_0}T(t,A)\in \mathcal{C}_p$ by hypothesis. For $n=1$, set $F(t)=\widetilde{B_0}T(t,A)\in \mathcal{C}_p$ and $G(t)=\tilde{B}T(t,A) \in \mathcal{C}_q$ which satisfy the hypotheses of Lemma~\ref{lem2}, respectively. Then, we get $H(t)=\widetilde{B_0}S_1(t)\in \mathcal{C}_{r_1}$ for $\frac{1}{p}+\frac{1}{q}=\frac{1}{r_1}$, which can be re-arranged into the expression of the theorem. For $n=2$, set $F(t)=\widetilde{B_0}T(t,A)\in \mathcal{C}_p$ and $G(t)=\tilde{B}S_{1}(t,A)\in\mathcal{C}_{q/2}$, to get $H(t)=\widetilde{B_0}S_2(t)\in \mathcal{C}_{r_2}$ where now
$
    \frac{2}{q}+\frac{1}{p}=\frac{1}{r_2}.
$ 
Continue this procedure until 
\[
    \max\left\{\frac{pq}{np+q},1\right\}=1.
\]
Clearing for the index, which swaps the value of the maximum, gives $n=\ell$.
After that, for $n>\ell$ we  can carry on applying Lemma~\ref{lem2} but now $r_n=1$. This ensures \ref{BtildeDyson_ii)}.

Now, let us show \ref{BtildeDyson_iii)}. Set $\psi_0(t)=\|\widetilde{B_0}S_n(t)\|_p$. Set $\psi_k(t)=\|\widetilde{B_0}S_{k-1}(t)\|_{\max\{1,\frac{q}{k}\}}$ as in Theorem~\ref{dyson2}. Then, all these real-valued functions are continuous and 
\[
\psi_n(t)\leq \psi_1*\psi_{n-1}(t)\leq \cdots \leq \psi_1^{[*n]}(t).
\]
For $n\in \mathbb{N}$, we have
\begin{equation} \label{boundwithtildeB0}
   \|\widetilde{B_0}S_n(t)\|_{\max\{1,r_n\}}\leq \psi_0*\psi_1^{[*n]}(t),
\end{equation}
in particular the index in the norm is $1$ for $n\geq \ell+1$.
Let
\[
   \theta_1(t)=\sum_{n=\ell+1}^\infty \psi_{1}^{[*n]}(t),
\]
which differs from $\theta_0(t)$ in the proof of Theorem~\ref{dyson2} by only a finite number of terms and so we know is convergent for all $t>0$ uniformly in $t\in (0,\alpha)$ for all $\alpha>0$. Then, from \eqref{boundwithtildeB0}, it follows that 
\begin{equation} \label{tailboundB0}
     \sum_{n=\ell+1}^\infty \|\widetilde{B_0}S_n(t)\|_1\leq
   \psi_0*\theta_1(t),
\end{equation}
so the series representation for $W(t)$ is indeed  convergent in $\mathcal{C}_1$ absolutely and uniformly in $t\in (0,\alpha)$. This ensures the validity of \ref{BtildeDyson_iii)}.

Now, consider \ref{BtildeDyson_iv)}. The case $n=0$ is an immediate consequence of the hypothesis. For the case $n=1,\ldots,\ell$, recall the argument invoking Lemma~\ref{lem2} in the proof of    \ref{BtildeDyson_ii)} above. From \eqref{expoestimate}, applied twice, once with $\psi=\psi_0$ and the other with $\psi=\psi_1$, we have
\[
     \int_0^\infty \mathrm{e}^{-\omega s}(\psi_0(s)+\psi_1(s))\,\mathrm{d}s<\infty,
\]       
whenever $\omega>\omega_0(A)$ is large enough. Hence, using \eqref{boundwithtildeB0}, 
\[
   \int_0^\infty \mathrm{e}^{-\omega s} \|\widetilde{B_0}S_n(s)\|_{r_n}\mathrm{d}s\leq
\left(\int_0^\infty \mathrm{e}^{-\omega s} \psi_0(s)\mathrm{d}s \right)
\left( \int_0^\infty \mathrm{e}^{-\omega s} \psi_1(s)\mathrm{d}s \right)^n<\infty.
\]
This ensures \ref{BtildeDyson_iv)}.

Finally, note that \ref{BtildeDyson_v)} follows from the fact that
\[
    \int_0^1 \psi_0*\theta_1(s)\mathrm{d}s<\infty,
\]
as a consequence of \cite[Lemma~13.4.3]{HP1957} and \eqref{tailboundB0}.
\end{proof}

\begin{Cor}
Let $T(t,A)$ be an immediately norm continuous one-parameter semigroup. For $p\geq q\geq 1$, let $B_0\in \mathcal{B}_p(A)$ and $B \in \mathcal{B}_{q}(A)$. Then, $B_0\in \mathcal{B}_{p}(A+B)$.
\end{Cor}
\begin{proof}
We verify Definition~\ref{d3} for index $p$ and generator $A+B$. Since $B_0,\,B\in \mathcal{B}(A)$, property \ref{a)forBq(A)} is straightforward. The property \ref{b)forBq(A)} follows from the expansion \eqref{BtildeDysonformula} and the fact that all the terms on the right hand side lie in $\mathcal{C}_p$. The property \ref{c)forBq(A)} is a consequence of \ref{BtildeDyson_iv)} and \ref{BtildeDyson_v)} in Theorem~\ref{BtildeDyson}.
\end{proof}

\begin{Cor}
The relation given in Definition~\ref{equivalentclass} is an equivalence. 
\end{Cor}
\begin{proof} 
Reflexivity follows from the fact that $0 \in \mathcal{B}_{q}(A)$ for any generator $A$. Next, consider the property of symmetry. 
Assume that $A_{1} \approx_q A_{2}$. This means that $B=A_{2}- A_{1} \in \mathcal{B}_{q}(A_{1})$. Clearly $-B \in  \mathcal{B}_{q}(A_{1})$. Moreover, since $-B\in \mathcal{J}(A_1)$, then $D(A_{1})=D(A_{2})$ and $-B \in  \mathcal{J}(A_{2})$. Now, according to Theorem~\ref{BtildeDyson}, taking $p=q$,  we have \[\int_{0}^{1} \|\widetilde{-B}T(s,A_{2})\|_{q} \textnormal{d}s<\infty.\] Thus $-B \in \mathcal{B}_{q}(A_{2})$ for $A_{1}=A_{2}+(-B)$ and  so indeed $A_{2} \approx_q A_{1}$. 

Finally, let us prove transitivity. Assume that $A_{1} \approx_q A_{2}$ and $A_{2} \approx_q A_{3}$. Then, $B_{1} = A_{2}-A_{1} \in \mathcal{B}_{q}(A_{1})$ and $B_{2} = A_{3}-A_{2} \in \mathcal{B}_{q}(A_{2})$. By symmetry, we now know that $\mathcal{B}_{q}(A_{2})=\mathcal{B}_{q}(A_{1})$. Therefore, because $\mathcal{B}_{q}(A_{1})$ is a linear space, $B_{1}+B_{2} \in \mathcal{B}_{q}(A_{1})$ too. Since $A_{3}=A_{1}+(B_{1}+B_{2})$, it follows that $A_{1} \approx_q A_{3}$ as needed.
\end{proof}

\section{The resolvent} \label{closed1}

In this section we examine how the resolvents of two generators which are $\mathcal{B}_q$-equivalent relate to one another. Our starting point is the fact that for $A_1\approx_q A_2$, the one-parameter semigroups $T(t,A_j)$ satisfy the spectral mapping theorem,
\[
     \mathrm{e}^{-t\operatorname{Spec}(A_j)}=\operatorname{Spec}(T(t,A_j))
\]
and
\[
     -\omega_0(A_j)=\inf \Re\left(\operatorname{Spec}(A_j)\right).
\]
Moreover
\begin{equation} \label{frominc}
     \lim_{|y|\to\infty}\|R(x+iy,A_j)\|_{\infty}=0 \qquad\forall x<-\omega_0(A_j).
\end{equation}
All this is a consequence of (and the latter is equivalent to) the fact that $T(t,A_j)$ are immediately norm continuous, see~\cite[Corollary~2.3.6]{vN1996} and \cite{Y1992}.  In fact the limit on the left hand side of \eqref{frominc} is zero for all $x\in \mathbb{R}$, see e.g. \cite[Theorem~3.6]{Pazy1983}.

Just as the class $\mathcal{B}(A)$ comprise relatively bounded perturbations of the generator $A$, we will see that the classes $\mathcal{B}_q(A)$ comprise relatively Schatten-class perturbations of $A$. Further, the norm of the composition of a $B\in \mathcal{B}_q(A)$ with the resolvent of $A$ goes to zero in lines parallel to the imaginary axis.

\begin{Lma} \label{relatpert}
Let $A$ be the generator of an immediately norm continuous semigroup and let $B\in \mathcal{B}_q(A)$. Then $\|B R(w,A)\|_q<\infty$ for all $w\not\in \operatorname{Spec}(A)$. Moreover,
\[
    \lim_{|y|\to \infty} \|B R(x+i y,A)\|_q=0 \qquad \forall x<-\omega_0(A).
\]
\end{Lma}
\begin{proof}
Let $z=x+iy$. We consider the proof of the first claim for $w=z$. If $x<-\omega_0(A)$, by \cite[Lemma~13.3.4]{HP1957},
\[
   BR(z,A)=\int_0^\infty \mathrm{e}^{zs}\tilde{B}T(s,A) \,\mathrm{d}s,
\]
where the integral converges in $\|\cdot\|_\infty$. Now for any $t>0$,
\begin{align*}
BR(z,A) & = \int_0^t \mathrm{e}^{zs}\tilde{B}T(s,A)\,\mathrm{d}s+
\mathrm{e}^{tz}\int_t^\infty \mathrm{e}^{z(s-t)}\tilde{B}T(s,A)\,\mathrm{d}s \\
&= \int_0^t \mathrm{e}^{zs}\tilde{B}T(s,A)\,\mathrm{d}s + \mathrm{e}^{tz} \int_0^\infty \mathrm{e}^{zs}\tilde{B}T(t+s,A)\,\mathrm{d}s \\
&= \int_0^t \mathrm{e}^{zs}\tilde{B}T(s,A)\,\mathrm{d}s+\mathrm{e}^{zt}\tilde{B}T(t,A)R(z,A).
\end{align*}
Hence,
\[
    \|BR(z,A)\|_q\leq \max\{1,\mathrm{e}^{xt}\} \left(\int_0^t\|\tilde{B}T(s,A)\|_q\,\mathrm{d}s + \|\tilde{B}T(t,A)\|_q \|R(z,A)\|_{\infty}\right)<\infty.
\]
Now let $w\not\in \operatorname{Spec}(A)$ such that $w\not=z$. As
\[
    BR(w,A)=BR(z,A)(I+(z-w)R(w,A)),
\]
the first claim follows.

For the second claim, we saw already that
\begin{align*}
\mathrm{e}^{-zt}BR(z,A) &= \mathrm{e}^{-zt}\int_0^t \mathrm{e}^{zs}\tilde{B}T(s,A)\,\mathrm{d}s+BR(z,A)T(t,A).
\end{align*}
Then 
\[
   \mathrm{e}^{-xt} \|BR(z,A)\|_{q}\leq \mathrm{e}^{-xt} \left\| \int_0^t \mathrm{e}^{zs}\tilde{B}T(s,A)\,\mathrm{d}s \right\|_q+\|BR(z,A)\|_q\|T(t,A)\|_{\infty}.  
\]
Assuming without loss of generality that $\omega_0(A)=0$ and letting
\[
    t>\frac{\log(M)}{-x}>0
\]
for $\|T(t,A)\|\leq M$, yields
\[
     \|BR(z,A)\|_q \leq \frac{\mathrm{e}^{-xt}}{\mathrm{e}^{-xt}-M}
\left\|\int_0^t \mathrm{e}^{iys} \mathrm{e}^{xs}\tilde{B}T(s,A)\,\mathrm{d}s\right\|_q.
\]
By virtue of Lemma~\ref{lemma1}, we know that
\[
   s\longmapsto \chi_{(0,t]}(s)\mathrm{e}^{xs}\tilde{B}T(s,A)
\]
is piecewise continuous in $\|\cdot\|_q$ and it lies in $L^1(\mathbb{R};\mathcal{C}_q)$ by the definition of $\mathcal{B}_q(A)$. Hence, the second conclusion follows from a version of the Riemann-Lebesgue lemma for $\mathcal{C}_q$-valued functions \cite[Theorem~C.8]{MR1721989}.
\end{proof}

From the first conclusion of the above lemma follows that, for $A_1\approx_q A_2$, the essential spectra of $A_1$ and $A_2$ coincide. In fact, since
\[
     R(z,A_2)-R(z,A_1)=R(z,A_2)BR(z,A_1),
\]
for $z\not\in \operatorname{Spec}(A_1)\cup \operatorname{Spec}(A_2)$, then
\begin{equation} \label{belowgrowthbound}
     \lim_{|y|\to \infty}\|R(x+iy,A_2)-R(x+iy,A_1)\|_q=0 \qquad \forall x<-\omega_0(A_j).
\end{equation}
Note that here there is no spectrum of $A_j$ for large enough $|y|$. Our next theorem strengthens these claims. 

\begin{Thm} \label{growthresolvents}
   Let $A_1$ and $A_2$ be two generators of immediately norm continuous semigroups, such that $A_1\approx_q A_2$.  
\begin{enumerate}[a)]
\item \label{growth1} Then, there exists a function $F:[0,\infty)\longrightarrow [0,\infty)$ such that
\[
    \operatorname{Spec}(A_1)\cup\operatorname{Spec}(A_2) \subset \{\lambda \in \mathbb{C}\,:\, |\operatorname{Im} \lambda |\leq F(|\operatorname{Re} \lambda|)\}.
\]
\item \label{growth2} Moreover, \eqref{belowgrowthbound} holds true for all $x\in \mathbb{R}$. 
\end{enumerate}
\end{Thm}
\begin{proof} 
We begin the proof by recalling a general property of the resolvent norm. Let $A$ be any closed operator on $\mathcal{H}$. If there exists $x\in \mathbb{R}$ such that 
\[
    \{y\in \mathbb{R}: x+iy\in \operatorname{Spec}(A)\}
\]
is a bounded set and
\begin{equation} \label{decayinsegment}
\lim_{|y|\to \infty} \|R(x+iy,A)\|_{\infty}=0,
\end{equation}
then there exists $F:[0,\infty)\longrightarrow [0,\infty)$ such that
\begin{equation} \label{universalenclosurespec}
    \operatorname{Spec}(A) \subset \{\lambda \in \mathbb{C}\,:\, |\operatorname{Im} \lambda |\leq F(|\operatorname{Re} \lambda|)\}
\end{equation}
and \eqref{decayinsegment} holds true for all $x\in \mathbb{R}$. 
As they are useful in their own right, we include here a self-contained proof of these statements. 

Since
\[
      \|R(x+iy,A)\|_{\infty} \geq \frac{1}{\operatorname{dist}\left(x+iy,\operatorname{Spec}A\right)},
\]
the hypothesis \eqref{decayinsegment} implies that 
\[
   \lim_{|y|\to \infty}\operatorname{dist}\left(x+iy,\operatorname{Spec}A\right)=\infty.
\]
Then, for any $w\in \mathbb{R}$, there exists $l(w)<\infty$ such that
\[
      \{y\in \mathbb{R}: w+iy\in \operatorname{Spec}(A)\}\subseteq \{w+iy\in\mathbb{C}\,:\, |y|\leq l(w)\},
\]
otherwise there is a contradiction to the above limit being infinity.
Setting $F(|w|)=\max\{l(w),\,l(-w)\}$ ensures \eqref{universalenclosurespec}.

Now let us show that \eqref{decayinsegment} holds for all $x\in \mathbb{R}$. Fix $w\in \mathbb{R}$. Then, \eqref{universalenclosurespec} ensures that the resolvent operator $R(w\pm iy,A)$ exists for $|y|$ large enough. Moreover,
\[
    R(w+iy,A)\left(1-(x-w)R(x+iy,A)\right)=R(x+iy,A)
\]
for such $y\in \mathbb{R}$. As $w$ is fixed, the hypothesis ensures that for $|y|$ large enough,
\[
     |x-w|\leq \frac{1}{2\|R(x+iy,A)\|}.
\]
Then
\[
    R(w+iy,A)=R(x+iy,A)\sum_{k=0}^\infty R(x+iy,A)^k(x-w)^k,
\]
where the right hand side is absolutely convergent in the operator norm. As,
\[
     \left\|\sum_{k=0}^\infty R(x+iy,A)^k(x-w)^k\right\|_{\infty}\leq 2
\]
for all such $y$, then the limit \eqref{decayinsegment} is also zero for $x=w$. 

We now complete the proof of the theorem. What we just showed ensures that \eqref{frominc} holds true for all $x\in \mathbb{R}$. Fix $x=\min\{-\omega_0(A_1),-\omega_0(A_2)\}-1$. Let $\lambda=x+iy$ and $\mu=w+iy$ for any $w\in \mathbb{R}$. Since
\begin{align*}
    \left(I+(\mu-\lambda)R(\lambda,A_2)\right)&\left(R(\mu,A_2)-R(\mu,A_1)\right) \\ &= 
\left(R(\lambda,A_2)-R(\lambda,A_1)\right)\left(1+(\lambda-\mu)R(\mu,A_1)\right),
\end{align*}
then 
\begin{align*} 
   R(\mu,A_2)&-R(\mu,A_1) \\ &=\left(1+(w-x)R(\lambda,A_2)\right)^{-1}\left(R(\lambda,A_2)-R(\lambda,A_1)\right)\left(1+(x-w)R(\mu,A_1)\right),
\end{align*}
 for $|y|$ large enough. By virtue of \eqref{frominc},
\begin{gather*}
    \lim_{|y|\to\infty} \left\| (1+(w-x)R(\lambda,A_2))^{-1}   \right\|_{\infty} =1 \qquad \text{and} \\
\lim_{|y|\to\infty} \left\| (1+(x-w)R(\mu,A_1))  \right\|_{\infty}=1.
\end{gather*}
Thus, \eqref{belowgrowthbound} also holds true for $x=w$.
\end{proof}

In the forthcoming sections, we consider applications of the above theory to perturbations of Schr\"odinger operators by a complex potential. 

\section{Perturbation of Schr\"odinger operators}

Let $-\Delta_{\Omega}$ denote the Dirichlet Laplacian on an open set $\Omega\subseteq \mathbb{R}^d$ for  $d=1,\,2,\,3$. Below we will consider two cases. Either $\Omega=\mathbb{R}^d$ or $\Omega$ is bounded, connected and its boundary is $C^{2}$. Set $\operatorname{Dom}(-\Delta_{\Omega})=W^{1,2}_{0}(\Omega)\cap W^{2,2}(\Omega)$. Then $-\Delta_{\Omega}$ is a self-adjoint non-negative operator and the generator of an immediately norm continuous semigroup. If $\Omega=\mathbb{R}^{d}$, we will simply write $-\Delta_{\mathbb{R}^d}=-\Delta$. 

Let $V\in L^{2}(\Omega)$ be a possibly complex-valued function. We will write the corresponding multiplication operator by $\mathrm{V}: \operatorname{Dom}(\mathrm{V}) \rightarrow L^{2}(\Omega)$ such that $\mathrm{V}:f \longmapsto Vf$, where $\operatorname{Dom}(\mathrm{V})=\{f \in L^{2}(\Omega) | \,Vf\in L^{2}(\Omega)\}$ is the maximal domain. Below we will use the same symbol $\mathrm{V}$ to denote this operator restricted to smaller domains.

Our main goal below will be to find conditions on $V$, so that $\mathrm{V}$ is a class $\mathcal{B}_{q}$ perturbation of $-\Delta_{\Omega}$ for some $1\leq q<\infty$. For this purpose, we begin by recalling the notion of the $l^{p}(L^{2}(\mathbb{R}^{d}))$ function spaces introduced by Birman and Solomjak, see \cite[Chapter~4]{MR2154153} for the original sources.

For $p\geq 1$ we say that a function $f\in L^{2}_{\mathrm{loc}}(\mathbb{R}^{d})$ is in $l^{p}(L^{2}(\mathbb{R}^{d}))$, if 
\begin{equation*}
\|f\|_{2;p}= \bigg(\sum_{\boldsymbol{\beta} \in \mathbb{Z}^{d}} \|\chi_{\boldsymbol{\beta}}f\|_{L^{2}(\mathbb{R}^{d})}^{p}\bigg)^{1/p}<\infty.
\end{equation*}
Here $\chi_{\boldsymbol{\beta}}$ is the characteristic function of the unit cube in $\mathbb{R}^{d}$ with center at $\boldsymbol{\beta}$. Note that
\begin{equation} \label{birmo1}
l^{p}(L^{2}(\mathbb{R}^{d}))\subseteq L^{2}(\mathbb{R}^{d}) \qquad \text{for all } \quad 1\leq p \leq 2.
\end{equation}
Moreover, for $\delta>0$, let $L^{2}_{\delta}(\mathbb{R}^{d})$ be the Lebesgue space of functions $f:\mathbb{R}^{d}\longrightarrow \mathbb{C}$ such that \[(1+|\cdot|^{2})^{\delta/2}f \in  L^{2}(\mathbb{R}^{d})\] with its standard norm. Then, \cite[(4.17)]{MR2154153} 
\begin{equation} \label{birmi}
L^{2}_{\delta}(\mathbb{R}^{d}) \subseteq l^{1}(L^{2}(\mathbb{R}^{d})) \qquad
\text{for} \qquad \delta>d/2.
\end{equation}
 
Let us describe the main tool in our arguments below from a general perspective. Let $H_{0}$ be the generator of an immediately norm continuous semigroup $\mathrm{e}^{-H_{0}t}=T(H_{0},t)$ on $L^{2}(\mathbb{R}^{d})$. Assume that $\mathrm{e}^{-H_{0}t}$ is associated to a heat kernel $K_{t}(\mathbf{x},\mathbf{y})$, so that for all $f\in L^{2}(\mathbb{R}^{d})$,
\begin{equation}\label{kerneldef}
\mathrm{e}^{-H_{0}t}f(\mathbf{x})=\int_{\mathbb{R}^{d}}K_{t}(\mathbf{x},\mathbf{y})f(\mathbf{y})\,\textnormal{d}\mathbf{y}.
\end{equation}
Further, assume that $K_{t}(\mathbf{x},\mathbf{y})$ satisfies the following Gaussian estimate for prescribed $b, k(t)>0$,
\begin{equation} \label{kernel}
|K_{t}(\mathbf{x},\mathbf{y})|\leq k(t) \mathrm{e}^{-b|\mathbf{x}-\mathbf{y}|^2/t} \qquad \forall \mathbf{x},\mathbf{y}\in \mathbb{R}^d, t>0,
\end{equation} 
where $k$ is such that 
\begin{equation}\label{kernel2}
     \int_{0}^{1} k(s)s^{\frac{d}{4}} \, \textnormal{d}s<\infty.
\end{equation}
Our goal is to show that  $\mathrm{V}$ is a $\mathcal{B}_{q}$ perturbation of $H_0$. 

\begin{Thm} \label{schroed}
Let $H_{0}$ be as above. If $V\in L^{2}(\mathbb{R}^{d})$, then $\mathrm{V}\in\mathcal{B}_{2}(H_{0})$.
\end{Thm}
\proof
We first show that $\mathrm{e}^{-H_{0}t}f\in \operatorname{Dom}(\mathrm{V})$ for all $f \in L^{2}(\mathbb{R}^{d})$. Indeed, using \eqref{kerneldef},
\begin{align*}
\int_{\mathbb{R}^{d}}|V(\mathbf{x})\mathrm{e}^{-H_{0}t}f(\mathbf{x})|^{2}\,\textnormal{d}\mathbf{x} =\int_{\mathbb{R}^{d}}|V(\mathbf{x})|^{2}\bigg|\int_{\mathbb{R}^{d}}K_{t}(\mathbf{x},\mathbf{y})f(\mathbf{y})\,\textnormal{d}\mathbf{y}\bigg|^{2}\,\textnormal{d}\mathbf{x}.
\end{align*} 
From the Cauchy--Schwarz inequality applied to the integral in $\mathbf{y}$, it follows that
\begin{align*}
\int_{\mathbb{R}^{d}}|V(\mathbf{x})\mathrm{e}^{-H_{0}t}f(\mathbf{x})|^{2}\,\textnormal{d}\mathbf{x} \leq\|f\|_{L^{2}(\mathbb{R}^{d})}^{2} \int_{\mathbb{R}^{d}}|V(\mathbf{x})|^{2}\int_{\mathbb{R}^{d}}|K_{t}(\mathbf{x},\mathbf{y})|^{2}\,\textnormal{d}\mathbf{y}\,\textnormal{d}\mathbf{x}.
\end{align*} 
Since the integrands are non-negative and measurable, Fubini's theorem applies, therefore we get
\begin{align}\label{lamosfa}
\int_{\mathbb{R}^{d}}|V(\mathbf{x})\mathrm{e}^{-H_{0}t}f(\mathbf{x})|^{2}\,\textnormal{d}\mathbf{x}\leq\|f\|_{L^{2}(\mathbb{R}^{d})}^{2} \int_{\mathbb{R}^{d}}\int_{\mathbb{R}^{d}}|V(\mathbf{x})K_{t}(\mathbf{x},\mathbf{y})|^{2}\,\textnormal{d}\mathbf{x}\,\textnormal{d}\mathbf{y}. 
\end{align} 
We now need to show that the integral in the right hand side of \eqref{lamosfa} is finite. Using \eqref{kernel} and Fubini's theorem again,
\begin{align*}
 \int_{\mathbb{R}^{d}}\int_{\mathbb{R}^{d}}|V(\mathbf{x})K_{t}(\mathbf{x},\mathbf{y})|^{2}\,\textnormal{d}\mathbf{x}\,\textnormal{d}\mathbf{y} \leq k^2(t)\int_{\mathbb{R}^{d}}|V(\mathbf{x})|^{2}\bigg(\int_{\mathbb{R}^{d}}\mathrm{e}^{-2b|\mathbf{x}-\mathbf{y}|^2/t}\,\textnormal{d}\mathbf{y}\bigg)\,\textnormal{d}\mathbf{x}.
\end{align*}
By a change of variables, $\mathbf{z}=\mathbf{x}-\mathbf{y}$, we obtain
\begin{align*}
\int_{\mathbb{R}^{d}}\mathrm{e}^{-2b|\mathbf{x}-\mathbf{y}|^2/t}\,\textnormal{d}\mathbf{y} =\int_{\mathbb{R}^{d}}\mathrm{e}^{-2b|\mathbf{z}|^2/t}\,\textnormal{d}\mathbf{z}=\left(\frac{t\pi}{2b}\right)^{d/2}.
\end{align*}
Therefore, we see that 
\begin{align}\label{momnet}
 \int_{\mathbb{R}^{d}}\int_{\mathbb{R}^{d}}|V(\mathbf{x})K_{t}(\mathbf{x},\mathbf{y})|^{2}\,\textnormal{d}\mathbf{x}\,\textnormal{d}\mathbf{y} \leq k^2(t)\|V\|_{L^{2}(\mathbb{R}^{d})}^{2}\left(\frac{t\pi}{2b}\right)^{d/2} <\infty.
\end{align}
Hence, $V\mathrm{e}^{-H_{0}t}f \in L^{2}(\mathbb{R}^{d})$.

Now, the estimates above also show that $\mathrm{V}\mathrm{e}^{-H_{0}t}\in\mathcal{C}_{2}(L^{2}(\mathbb{R}^{d}))$ for all $t>0$. Since $\mathrm{e}^{-H_{0}t}$ is immediately norm continuous, we have that $\mathrm{V}\mathrm{e}^{-H_{0}t}$ is continuous in the norm of $\mathcal{C}_{2}(L^{2}(\mathbb{R}^{d}))$. Hence, $\|\mathrm{V}\mathrm{e}^{-H_{0}t}\|_{2}$ is strongly measurable and by \eqref{momnet} and \eqref{kernel2},
\begin{equation*}
\int_{0}^{1}\|\mathrm{V}\mathrm{e}^{-H_{0}s}\|_{2}\,\textnormal{d}s \leq \left(\frac{\pi}{2b}\right)^{d/4}\|V\|_{L^{2}(\mathbb{R}^{d})}\int_{0}^{1}k(s)s^{\frac{d}{4}}\,\textnormal{d}s <\infty.
\end{equation*}
Thus, \begin{equation}\label{serorm}\int_{0}^{1}\|\mathrm{V}\mathrm{e}^{-H_{0}s}\|\,\textnormal{d}s <\infty,\end{equation} so that $\operatorname{Dom}(H_{0})\subseteq \operatorname{Dom}(V)$, by \cite[Lemma~11.4.4]{MR2359869}, and $\mathrm{V}\in \mathcal{B}(H_0)$. 
Additionally, we have seen that $\mathrm{V}\mathrm{e}^{-H_{0}t}\in \mathcal{C}_{2}(L^{2}(\mathbb{R}^{d}))$ and $\int_{0}^{1}\|\mathrm{V}\mathrm{e}^{-H_{0}s}\|_{2}\,\textnormal{d}s<\infty$. Therefore, $\mathrm{V}\in \mathcal{B}_{2}(H_0)$ as claimed.
\endproof

The Laplacian and the wide variety of operators considered in \cite{MR670130} satisfy the condition \eqref{kernel} with
\[
    k(t)=C\mathrm{e}^{at}t^{-d/2},
\]  
for some $C>0$, $a\in \mathbb{R}$ which depend on the specific operator. See \cite[Proposition~B.6.7]{MR670130}. We highlight that condition \eqref{kernel2} is guaranteed only for $d=1,\,2,\,3$. 

\begin{Rem} We strongly suspect that the above method applies to perturbations of certain magnetic Schr\"odinger operators, as the integral kernels satisfy an estimate of the form \eqref{kernel}. This relies on the so-called diamagnetic inequality. We can use the terminology and results of \cite{MR1756112} and \cite{MR670130} to sketch this. Indeed, consider the operator $H_{0}(A,V)=\frac{1}{2}(-\mathrm{i}\nabla-A)^{2}+V$ where $A \in \mathcal{H}_{\mathrm{loc}}(\mathbb{R}^{d})$, $V \in \mathcal{K}_{\pm}(\mathbb{R}^{d})$ as defined in \cite{MR1756112}. The heat kernel has the property
\begin{equation*}
|k_{t,A}(x,y)| \leq k_{t,0}(x,y) \qquad \forall x,y \in \mathbb{R}^{d} \text{ and } t>0.
\end{equation*}
See \cite[(6.30)]{MR1756112}. This, along with \cite[Theorem~B.7.1]{MR670130}, which expresses that $ k_{t,0}(x,y)$ satisfies an estimate of the form \eqref{kernel}, imply that $k_{t,A}(x,y)$ does as well.
\end{Rem}

\section{Characterisation of $\mathcal{B}_{q}$ perturbations of the Laplacian on $\mathbb{R}^d$}
If we consider the particular case of the Laplacian on $L^2(\mathbb{R}^{d})$, we can use certain results obtained via the functional calculus for self-adjoint operators to gain additional information on its $\mathcal{B}_{q}$ perturbations. In particular, we will use some of the statements described in \cite[Chapter~4]{MR2154153}, which we now recall, along with details of their applicability to the current perturbation theory.

Formally, for two functions $V,g:\mathbb{R}^{d}\rightarrow \mathbb{R}$ in some adequate function space, we can consider the operator $\mathrm{V}\mathrm{g}(-\mathrm{i}\nabla)$ on $L^2(\mathbb{R}^{d})$. Here $\mathrm{g}(-\mathrm{i}\nabla)$ acts as \[\mathrm{g}(-\mathrm{i}\nabla)f=\mathcal{F}^{-1}(g\mathcal{F}(f)),\] where $\mathcal{F}, \mathcal{F}^{-1}$ are the Fourier and inverse Fourier transforms, see \cite[Section~3.1]{MR2359869}.  We recall the following.
\begin{enumerate}[a)]
\item \label{boring1} If $V,g$ are both non-zero and $\mathrm{V}\mathrm{g}(-\mathrm{i}\nabla)\in \mathcal{C}_{2}(L^2(\mathbb{R}^{d}))$, then $V,g \in L^{2}(\mathbb{R}^{d})$. See \cite[Proposition~4.4]{MR2154153}.
\item \label{boring2} Let $1\leq p\leq 2$. If $V,g \in l^{p}(L^{2}(\mathbb{R}^{d}))$, then $\mathrm{V}\mathrm{g}(-\mathrm{i}\nabla)\in \mathcal{C}_{p}(L^2(\mathbb{R}^{d}))$ and in particular 
\begin{equation}\label{boring5}
\|\mathrm{V}\mathrm{g}(-\mathrm{i}\nabla)\|_{p}\leq C_{p}\|V\|_{2;p}\|g\|_{2;p},
\end{equation}
for some $C_{p}>0$. See \cite[Theorem~4.5]{MR2154153}.
\item \label{boring3}  If $V,g$ are both non-zero and $\mathrm{V}\mathrm{g}(-\mathrm{i}\nabla)\in \mathcal{C}_{1}(L^2(\mathbb{R}^{d}))$, then $V,g \in l^{1}(L^{2}(\mathbb{R}^{d}))$. See \cite[Proposition~4.7]{MR2154153}.
\end{enumerate}

For $g(\mathbf{x})= \mathrm{e}^{-|\mathbf{x}|^{2}t}$ we know that $\mathrm{g}(-\mathrm{i}\nabla)=\mathrm{e}^{\Delta t}$. Therefore, the point \ref{boring2} above implies that, if $V \in l^{p}(L^{2}(\mathbb{R}^{d}))$ for some $1\leq p \leq 2$, then the operator $\mathrm{V}\mathrm{e}^{\Delta t}$ is in $\mathcal{C}_{p}(L^2(\mathbb{R}^{d}))$, with an explicit bound on its $\mathcal{C}_{p}(L^2(\mathbb{R}^{d}))$ norm depending on $t$. This yields that $\mathrm{V}\in \mathcal{B}_{p}(-\Delta)$. Here, the possible $p$ will be found to depend on the dimension $d$. The reason for this is that the $l^{p}(L^{2}(\mathbb{R}^{d}))$ norm of $g$ depends on $d$. Accordingly, the $t$-dependence of the bound for $\mathrm{V}\mathrm{e}^{\Delta t}$ that can be achieved by \eqref{boring5} will depend on $d$, affecting the integrability of $\|\mathrm{V}\mathrm{e}^{\Delta t}\|_{p}$.

In addition, we can use \ref{boring1} and \ref{boring3} above, to find exact characterisations of the $\mathcal{B}_{1}(-\Delta)$ and $\mathcal{B}_{2}(-\Delta)$ classes of perturbations. The following is a stronger result than Theorem~\ref{schroed} for the Laplacian. Its proof is essentially an application of \eqref{boring5}.

\begin{Thm} \label{simonresult1}
For $d\leq3$, consider the operator $-\Delta$ on $L^2(\mathbb{R}^{d})$ with $\operatorname{Dom}(-\Delta)=W^{2,2}(\mathbb{R}^{d})$. If $V\in l^{p}(L^{2}(\mathbb{R}^{d}))$ for $1\leq p \leq 2$ such that $p>d/2$, then $\mathrm{V}\in \mathcal{B}_{p}(-\Delta)$.
\end{Thm}
\proof 
By virtue of \eqref{birmo1}, following the proof of Theorem~\ref{schroed}, we gather that $\mathrm{V}\in \mathcal{B}(-\Delta)$. Now, for $t>0$, choosing $g=\mathrm{e}^{-|\cdot|^{2}t}$ in \eqref{boring5}, we see that for some $C_{p}>0$,
\begin{equation*}
\|\mathrm{V}\mathrm{e}^{\Delta t}\|_{p}\leq C_{p}\|V\|_{2;p} \|\mathrm{e}^{-|\cdot|^{2}t}\|_{2;p}.
\end{equation*}
Further, for any $p\geq 1$ and $t>0$, we have that $\mathrm{e}^{-|\cdot|^{2}t}\in l^{p}(L^{2}(\mathbb{R}^{d}))$ with
\begin{align*}
\|\mathrm{e}^{-|\cdot|^{2}t}\|_{2;p}< 2^{d/p}\sqrt{\frac{\pi}{p}} (1+t^{-1/2})^{d/p}.
\end{align*}
We omit the details of this calculation. They can be found in \cite[Lemma~5.3.1]{DimoudisPhDThesis2022}.  Then,
\begin{equation*}
\|\mathrm{V}\mathrm{e}^{\Delta t}\|_{p}< \tilde{C}_{p}\|V\|_{2;p}(1+t^{-1/2})^{d/p}.
\end{equation*}
In particular, the right hand side of this inequality is finite for $t>0$ and thus $\mathrm{V}\mathrm{e}^{\Delta t}\in \mathcal{C}_{p}(\mathcal{H})$ for $1\leq p\leq 2$. 
Finally, we have
\begin{equation*}
\int_{0}^{1}\|\mathrm{V}\mathrm{e}^{-H_{0}s}\|_{p} \,\textnormal{d}s< \tilde{C}_{p}\|V\|_{2;p}\int_{0}^{1}(1+s^{-1/2})^{d/p}  \,\textnormal{d}s. 
\end{equation*}
Noting that $1<s^{-1/2}$ for $s \in (0,1)$, then gives
\begin{equation*}
\int_{0}^{1}\|\mathrm{V}\mathrm{e}^{-H_{0}s}\|_{p} \,\textnormal{d}s<2^{d/p}\tilde{C}_{p}\|V\|_{2;p}\int_{0}^{1}s^{-d/2p}\,  \textnormal{d}s<\infty,
\end{equation*}
for $1\leq p\leq 2$ with $p>d/2$. Putting these observations together, we get that $\mathrm{V}\in \mathcal{B}_{p}(-\Delta)$.
\endproof

In the above theorem, $p\in[1,2]$ for $d=1$, but $p\in(1,2]$ for $d=2$. The following corollary is a consequence of \eqref{birmi}.

\begin{Cor} \label{simonresult2}
Consider the operator $-\frac{d^2}{dx^2}$ on $L^2(\mathbb{R})$ with $\operatorname{Dom}\left(-\frac{d^2}{dx^2}\right)=W^{2,2}(\mathbb{R})$. If $V \in L^{2}_{\delta}(\mathbb{R})$ for some $\delta>1/2$, then $\mathrm{V}\in \mathcal{B}_{1}\left(-\frac{d^2}{dx^2}\right)$.
\end{Cor}

We also highlight the following necessary and sufficient condition for multiplication operators to lie in $\mathcal{B}_{1}(-\Delta)$ or $\mathcal{B}_{2}(-\Delta)$.

\begin{Thm} \label{simonresult3}
Let $V:\mathbb{R}^{d}\rightarrow \mathbb{C}$ be non-zero.
\begin{enumerate}[i)]
\item \label{1vr} For $d\leq 3$, $\mathrm{V}\in \mathcal{B}_{2}(-\Delta)$ if and only if $V \in L^{2}(\mathbb{R}^{d})$.
\item \label{2vr} For $d=1$, $\mathrm{V}\in \mathcal{B}_{1}\left(-\frac{d^2}{dx^2}\right)$ if and only if $V \in l^{1}(L^{2}(\mathbb{R}))$.
\end{enumerate}
\end{Thm}
\proof
The forward directions of these two claims can be shown by using \ref{boring1} and \ref{boring3} above, respectively. These results imply that the condition~\ref{b)forBq(A)} in Definition~\ref{d3} is fulfilled, only if $V$ is in the respective function spaces. We show this for claim \ref{1vr}, the case of claim \ref{2vr} being analogous. If $\mathrm{V}\in \mathcal{B}_{2}(-\Delta)$, then $\mathrm{V}\mathrm{e}^{\Delta t} \in \mathcal{C}_{2}(L^2(\mathbb{R}^d))$. Therefore, \ref{boring1} above, implies that $V \in L^{2}(\mathbb{R}^{d})$.

The other directions of these two claims follow from theorems~\ref{schroed} and~\ref{simonresult1} respectively. 
\endproof

\section{Dirichlet Laplacian on a bounded region}

In this final section, we consider another application of the semigroup theory we have described in this work, in the spirit of the results developed in \cite{B2019}. We will derive eigenvalue asymptotics for non-self-adjoint perturbations of the Dirichlet Laplacian on $\Omega$ bounded, open and connected. In order to simplify technical details, we assume additionally that the boundary of $\Omega$ is $C^2$. Below we denote the Lebesgue measure of $\Omega$ by $|\Omega|$.

The Dirichlet Laplacian on $\Omega$ is the generator of a Gibbs semigroup. More specifically, $\mathrm{e}^{\Delta_{\Omega}t}$ has a positive integral kernel $K_{\Omega,t}(\mathbf{x},\mathbf{y})$ which is in $C^{\infty}(\mathbb{R}^{d}\times \mathbb{R}^{d})$ and which satisfies the following Gaussian estimate
\begin{equation} \label{gaussianfinite}
K_{\Omega,t}(\mathbf{x},\mathbf{y})\leq \frac{1}{(4\pi t)^{d/2}} \mathrm{e}^{-|\mathbf{x}-\mathbf{y}|^2/4t} \qquad \forall \mathbf{x},\mathbf{y}\in \Omega, \,t>0.
\end{equation}
As $-\Delta_{\Omega}$ is self-adjoint and has a compact resolvent, its spectrum is purely discrete. We write the eigenvalues of $-\Delta_{\Omega}$ in non-decreasing order as $\{\mu_{n}\}_{n=1}^{\infty}$. Recall that \cite[Theorem~6.3.1]{MR1349825}, there exists a constant $a_0(\Omega)>0$ such that 
\begin{equation*} 
\mu_{n}\geq n^{2/d} a_0(\Omega) \qquad \forall n\in \mathbb{N}.
\end{equation*}
We now show how to replicate this estimate asymptotically, for the real part of the eigenvalues of  $-\Delta_{\Omega}+\mathrm{V}$, where $V\in L^{2}(\Omega)$.
The proof of this is similar to that of \cite[Corollary~3]{B2019}. The fact that the spectrum of the perturbed operator is countably infinite is part of the conclusion. 

\begin{Thm}\label{somethingsome}
Let $d\leq 3$. Let $\Omega\subseteq \mathbb{R}^d$ be bounded, open and connected. Further, let $V\in L^{2}(\Omega)$. Then, there exists an infinite sequence $\{\lambda_{k}\}_{k=1}^{\infty}$ such that $\sigma(-\Delta_{\Omega}+\mathrm{V})=\{\lambda_{k}\}_{k=1}^{\infty}$. Moreover, there exists an $N\in \mathbb{N}$, such that for $n\geq N$, 
\begin{equation} \label{sub-Weyl}
\Re(\lambda_{n})\geq \frac{4\pi}{(4\mathrm{e}|\Omega|)^{2/d}} n^{2/d}.\end{equation}
\end{Thm}
\proof
We aim to prove that $\mathrm{V}\in \mathcal{B}_2(-\Delta_{\Omega})$, and then proceed as in the proof of \cite[Corollary~3]{B2019}. 

Let $f\in L^{2}(\Omega)$. Then, similarly to the proof of \eqref{lamosfa}, for $t>0$
\begin{align*}
\int_{\Omega}|V(\mathbf{x})\mathrm{e}^{\Delta_{\Omega}t}f(\mathbf{x})|^{2}\,\textnormal{d}\mathbf{x}& =\int_{\Omega}|V(\mathbf{x})|^{2}\bigg|\int_{\Omega}K_{\Omega,t}(\mathbf{x},\mathbf{y})f(\mathbf{y})\,\textnormal{d}\mathbf{y}\bigg|^{2}\,\textnormal{d}\mathbf{x} \\ & \leq\|f\|_{L^{2}(\Omega)}^{2} \int_{\Omega}|V(\mathbf{x})|^{2}\int_{\Omega}|K_{\Omega,t}(\mathbf{x},\mathbf{y})|^{2}\,\textnormal{d}\mathbf{y}\,\textnormal{d}\mathbf{x} \\&=\|f\|_{L^{2}(\Omega)}^{2} \int_{\Omega}\int_{\Omega}|V(\mathbf{x})K_{\Omega,t}(\mathbf{x},\mathbf{y})|^{2}\,\textnormal{d}\mathbf{x}\,\textnormal{d}\mathbf{y}. 
\end{align*} 
We now prove that the integral $\int_{\Omega}\int_{\Omega}|V(\mathbf{x})K_{\Omega,t}(\mathbf{x},\mathbf{y})|^{2}\,\textnormal{d}\mathbf{x}\,\textnormal{d}\mathbf{y}$ is convergent. Indeed, from the Gaussian estimate \eqref{gaussianfinite}, we see that
\begin{equation*}\label{approx1}
\int_{\Omega}\int_{\Omega}|V(\mathbf{x})K_{\Omega,t}(\mathbf{x},\mathbf{y})|^{2}\,\textnormal{d}\mathbf{x}\,\textnormal{d}\mathbf{y} \leq \frac{1}{(4\pi t)^{d}}\int_{\Omega}|V(\mathbf{x})|^{2}\bigg(\int_{\Omega}\mathrm{e}^{-|\mathbf{x}-\mathbf{y}|^2/2t}\,\textnormal{d}\mathbf{y}\bigg)\,\textnormal{d}\mathbf{x}.
\end{equation*}
Noting that \[\int_{\Omega}\mathrm{e}^{-|\mathbf{x}-\mathbf{y}|^2/2t}\,\textnormal{d}\mathbf{y}\leq \int_{\mathbb{R}^{d}}\mathrm{e}^{-|\mathbf{x}-\mathbf{y}|^2/2t}\,\textnormal{d}\mathbf{y},\] and following the proof of Theorem~\ref{schroed},
\begin{align*}\label{approx1}
\frac{1}{(4\pi t)^{d}}\int_{\Omega}|V(\mathbf{x})|^{2}\bigg(\int_{\Omega}\mathrm{e}^{-|\mathbf{x}-\mathbf{y}|^2/2t}\,\textnormal{d}\mathbf{y}\bigg)\,\textnormal{d}\mathbf{x} & \nonumber \leq \frac{1}{(4\pi t)^{d}} \|V\|_{L^{2}(\Omega)}^{2}\left(\int_{\mathbb{R}^{d}}\mathrm{e}^{-|\mathbf{z}|^2/2t}\,\textnormal{d}\mathbf{z}\right)
\\&\leq \frac{1}{(4\pi t)^{d}}\|V\|_{L^{2}(\Omega)}^{2}\left(2\pi t\right)^{d/2} \\ \nonumber&=\frac{1}{(8\pi t)^{d/2}}\|V\|_{L^{2}(\Omega)}^{2}< \infty.
\end{align*}
These calculations show that $\mathrm{V}\mathrm{e}^{\Delta_{\Omega}t}f\in L^{2}(\Omega)$. In other words, 
\[\bigcup_{t>0} \mathrm{e}^{\Delta_{\Omega}t}(L^{2}(\Omega)) \subseteq\operatorname{Dom}(\mathrm{V}).\]

Since the operator $\mathrm{V}\mathrm{e}^{\Delta_{\Omega}t}$ has integral kernel $V(\mathbf{x})K_{\Omega,t}(\mathbf{x},\mathbf{y})$, and we have seen that 
\begin{equation*}
\int_{\Omega}\int_{\Omega}|V(\mathbf{x})K_{\Omega,t}(\mathbf{x},\mathbf{y})|^{2}\,\textnormal{d}\mathbf{x}\,\textnormal{d}\mathbf{y}<\infty,
\end{equation*}
then $\|\mathrm{V}\mathrm{e}^{\Delta_{\Omega}t}\|_{2}<\infty$.
In addition, justifying the strong measurability of $\mathrm{V}\mathrm{e}^{\Delta_{\Omega}t}$ as in the proof of Theorem~\ref{schroed}, we gather that 
\begin{equation*}
\int_{0}^{1}\|\mathrm{V}\mathrm{e}^{\Delta_{\Omega}s}\|_{\infty}\,\textnormal{d}s\leq \int_{0}^{1}\|\mathrm{V}\mathrm{e}^{\Delta_{\Omega}s}\|_{2}\,\textnormal{d}s \leq \frac{\|V\|_{L^{2}(\Omega)}}{(8\pi)^{d/4}}\int_{0}^{1}s^{-d/4}\,\textnormal{d}s <\infty,
\end{equation*}
for $d\leq 3$.
Therefore, $\operatorname{Dom}(-\Delta_{\Omega})\subseteq \operatorname{Dom}(\mathrm{V})$ and $\mathrm{V}\in \mathcal{B}_{2}(-\Delta_{\Omega})$. 

Now, by virtue of Corollary~\ref{gibbs}, it follows that $\mathrm{e}^{(\Delta_{\Omega}-\mathrm{V})t}$ is a Gibbs semigroup. Moreover, the triangle inequality alongside with Corollary~\ref{Duhamel}, imply that there exist $M_{1},\gamma>0$ such that
\begin{equation}\label{approx2} 
\|\mathrm{e}^{(\Delta_{\Omega}-\mathrm{V})t}\|_{2}\leq \|\mathrm{e}^{\Delta_{\Omega}t}\|_{2}+M_{1},
\end{equation}
for $0<t\leq \gamma$.

Since $\mathrm{V}$ is relatively compact, hence relatively bounded with bound zero with respect to $-\Delta_{\Omega}$, there exists an infinite sequence $\{\lambda_{k}\}_{k=1}^{\infty}\subset \mathbb{C}$ such that 
\[\operatorname{Spec}(-\Delta_{\Omega}+V)= \{\lambda_{k}\}_{k=1}^{\infty}.\] See for example \cite[Corollary~4.10]{MR3588930}. Moreover, \[\operatorname{Spec}(\mathrm{e}^{\left(\Delta_{\Omega}-\mathrm{V}\right)t})= \{0\}\cup\{\mathrm{e}^{-\lambda_{k}t}\}_{k=1}^{\infty},\] where 
\[
\lim_{n\rightarrow \infty}\Re(\lambda_{n})=\infty.
\] 
From this it follows that the sequence $\{\lambda_{k}\}_{k=1}^{\infty}$ can be reordered so that  $\Re(\lambda_{n})$ is non-decreasing. We assume that the latter is the case.
Then, 
\begin{align*}
\sum_{k=1}^{\infty}\mathrm{e}^{-\Re(\lambda_{k})t}&=\sum_{k=1}^{\infty}|\mathrm{e}^{-\lambda_{k}t}| \\ & \leq \|\mathrm{e}^{(\Delta_{\Omega}-\mathrm{V})t}\|_{1} \\ & =\|\mathrm{e}^{(\Delta_{\Omega}-\mathrm{V})t/2}\mathrm{e}^{(\Delta_{\Omega}-\mathrm{V})t/2}\|_{1} \\ &\leq \|\mathrm{e}^{(\Delta_{\Omega}-\mathrm{V})t/2}\|^{2}_{2} \\ &\leq 2\|\mathrm{e}^{\Delta_{\Omega}t/2}\|^{2}_{2}+2M_{1}^{2},
\end{align*}
for $0<t\leq \gamma$.  Here, we have used \eqref{approx2} and the inequality $(a+b)^{2}\leq 2a^{2}+2b^{2}$ for $a, b \in \mathbb{R}$ in the last step.

We now compute a bound for $\|\mathrm{e}^{\Delta_{\Omega}t/2}\|^{2}_{2}$, using \eqref{gaussianfinite} and the Fubini theorem. We have that
\begin{align}\label{approx3}
\nonumber \|\mathrm{e}^{\Delta_{\Omega}t/2}\|^{2}_{2}&=\int_{\Omega}\int_{\Omega}|K_{\Omega,t/2}(\mathbf{x},\mathbf{y})|^{2}\,\textnormal{d}\mathbf{x}\,\textnormal{d}\mathbf{y} \\ \nonumber&\leq \frac{1}{(2\pi t)^{d}}\int_{\Omega}\bigg(\int_{\Omega}\mathrm{e}^{-|\mathbf{x}-\mathbf{y}|^2/t}\,\textnormal{d}\mathbf{y}\bigg)\,\textnormal{d}\mathbf{x} \\ \nonumber&\leq \frac{1}{(2\pi t)^{d}} \int_{\Omega}\bigg(\int_{\mathbb{R}^{d}}\mathrm{e}^{-|\mathbf{z}|^2/t}\,\textnormal{d}\mathbf{z}\bigg)\,\textnormal{d}\mathbf{x}
\\ &=\frac{|\Omega|}{(4\pi t)^{d/2}}.
\end{align}
Hence,
\begin{equation}\label{weylll} 
\sum_{k=1}^{\infty}\mathrm{e}^{-\Re(\lambda_{k})t}\leq \frac{2|\Omega|}{(4\pi t)^{d/2}}+2M_{1}^{2},
\end{equation}
for $0<t\leq \gamma$. 

Now, since $\{\Re(\lambda_{k})\}_{k=1}^{\infty}$ is non-decreasing, 
\begin{equation*}
n\mathrm{e}^{-\Re(\lambda_{n})t}=\sum_{k=1}^{n}\mathrm{e}^{-\Re(\lambda_{n})t}\leq \sum_{k=1}^{n}\mathrm{e}^{-\Re(\lambda_{k})t}\leq \sum_{k=1}^{\infty}\mathrm{e}^{-\Re(\lambda_{k})t} \qquad \forall n \in \mathbb{N}.
\end{equation*}
Consequently,
\begin{equation}\label{jabjab}
n\mathrm{e}^{-\Re(\lambda_{n})t}\leq \frac{2|\Omega|}{(4\pi t)^{d/2}}+2M_{1}^{2}.
\end{equation}
Since $\lim_{t\rightarrow 0^{+}}t^{-d/2}=\infty$, there exists $\gamma_{1}\leq\gamma$ such that  $2M_{1}^{2}\leq \frac{2|\Omega|}{(4\pi t)^{d/2}}$ for all $0<t\leq \gamma_{1}$. Thus, 
\begin{equation}\label{jabojab}
n\mathrm{e}^{-\Re(\lambda_{n})t}\leq \frac{4|\Omega|}{(4\pi t)^{d/2}},
\end{equation}
for $0<t\leq \gamma_{1}$. Also, there exists $N\in \mathbb{N}$ such that for $n\geq N$,
\[\frac{1}{\Re(\lambda_{n})}<\gamma_{1}.\]
So, put $t=\frac{1}{\Re(\lambda_{n})}$ in \eqref{jabojab}. Then, for $n\geq N$,
\begin{equation*}
n\mathrm{e}^{-1}\leq \frac{4|\Omega|}{(4\pi)^{d/2}}\Re(\lambda_{n})^{d/2},
\end{equation*}
or 
\begin{equation*}
\Re(\lambda_{n})\geq \frac{4\pi}{(4\mathrm{e}|\Omega|)^{2/d}} n^{2/d},
\end{equation*}
which completes the proof.
\endproof

We mention that it is likely that, by following the line of arguments in the proof of \cite[Corollary~3]{B2019}, involving the holomorphic semigroup $\mathrm{e}^{\Delta_{\Omega}\tau}$ for $\tau$ in a half plane, one may be able to determine in a similar way as above, some  information about the asymptotic behaviour of the imaginary part of the eigenvalues of $-\Delta_{\Omega}+\mathrm{V}$.  

\begin{Rem}
It is not our aim above to obtain any optimal constant on the right hand side of \eqref{sub-Weyl}, but rather illustrate a perturbation method for the spectrum which is based on the theory of one-parameter semigroups. The classical Weyl asymptotic formulas, of which a significant amount of detail is known in the self-adjoint setting, see \cite{SaVa} and references therein, predict that the term $(4\mathrm{e})^{2/d}>10^{2/d}$ in the denominator could be improved, but should not be smaller than 
\[
\Gamma\left(1+\frac{d}{2}\right)^{2/d}=\begin{cases} \frac{\pi}{4} & d=1 \\ 1 & d=2 \\
\frac{3^{2/3}\pi^{1/3}}{2^{4/3}} & d=3 .\end{cases}
\]
However, note that the other terms and powers match the optimal coefficient of the classical self-adjoint case.
\end{Rem}

\section*{Acknowledgements}
SD was supported by Heriot-Watt University's James Watt scholarships programme.

\printbibliography

\end{document}